\documentclass[11pt]{amsart}
\usepackage{graphicx}
\usepackage{amsmath,amsfonts,amssymb,latexsym,amsthm,amscd}
 \newtheorem{theorem}{Theorem}[section]

\begin{document}

\title[Dehn Twist Subgroup Of Non-orientable Surface]
{The Torsion Generating Set Of The Dehn Twist Subgroups Of Non-orientable Surfaces}

\author{Xiaoming Du}
\address{South China University of Technology,
  Guangzhou 510640, P.R.China}
\email{scxmdu@scut.edu.cn}



\thanks{The author would like to thank Szepietowski for telling to him the generators of the index 2 Dehn twist subgroup
and pointing out some generator was missing in the earlier version of the paper.}

\maketitle

\begin{abstract}
Let $N_g$ be the non-orientable surface with genus $g$, $\text{MCG}(N_g)$ be the mapping class group of $N_g$,
$\mathcal{T}(N_g)$ be the index 2 subgroup generated by all Dehn twists of $\text{MCG}(N_g)$.
We prove that for odd genus, $\mathcal{T}(N_g)$ can be generated by three elements of finite orders.
\end{abstract}

\section{Introduction}

Let $N_g$ be the non-orientable surface with genus $g$, $\text{MCG}(N_g)$ be the mapping class group of $N_g$.
Lickorish was the first one to discover that all Dehn twists can only generate an index 2 subgroup of $\text{MCG}(N_g)$ (\cite{Li}).
We denote this subgroup as $\mathcal{T}(N_g)$. Outside $\mathcal{T}(N_g)$,
there is a mapping class called "Y-homeomorphism" or "cross-cap slide".
A finite set of generators for $\text{MCG}(N_g)$ and $\mathcal{T}(N_g)$ was given by Chillingworth (\cite{Ch}).
When $g=2$, Lickorish found $\text{MCG}(N_2) \cong \mathbb{Z}_2 \oplus \mathbb{Z}_2$,
Chillingworth found $\mathcal{T}(N_2)$ can be generated by one Dehn twist.
When $g=3$, Birman and Chillingworth proved that $\text{MCG}(N_3)$ can be generated by three involutions (\cite{BC}),
Chillingworth found $\mathcal{T}(N_3)$ can be generated by two Dehn twists.

It is a natural question how to simplify the generating sets
for $\text{MCG}(N_g)$ and $\mathcal{T}(N_g)$ as much as possible when $g$ is large.
We want to reduce both the number and the orders of the generators.
When $g \geq 4$, a generating set for $\text{MCG}(N_g)$ consisting of four involutions was constructed by Szepietowski.
Szepietowski also proved when $g \geq 4$, $\text{MCG}(N_g)$ can be generated by three elements (See \cite{Sz}).
The first homology of $\text{MCG}(N_g)$ has been calculated by Korkmaz.
By Korkmaz's result, when $g=4$, the smallest number of generators for $\text{MCG}(N_4)$ is at least 3,
So the minimal number of the generators for $\text{MCG}(N_3)$ is 3.
About $\mathcal{T}(N_g)$,
Stukow gave a finite presentation of $\mathcal{T}(N_g)$ in \cite{St}.
Omori reduced the number of Dehn twist generators for $\mathcal{T}(N_g)$ to $g+1$ when $g \geq 4$ (\cite{Om}).

In \cite{Du1}, the author proved the following: when the genus $g' \geq 5$,
the extended mapping class group $\text{MCG}^{\pm}(S_{g'})$ can be generated by two elements of finite order.
One is of order 2 and the other is of order $4g'+2$.
In \cite{Du2} (preprint), the author proved that the above result is also true for $g'=3,4$.
We found that the method in \cite{Du1} \cite{Du2} can be used in the case of $\mathcal{T}(N_g)$.
We have the following result:

\begin{theorem}
Let $\mathcal{T}(N_g)$ be the index 2 Dehn twist subgroup of the mapping class group of a non-orientable surface.
If $g \geq 5$ is odd, $\mathcal{T}(N_g)$ can be generated by three elements of finite order.
One of the generator is of order $2g$. The other two are of order $2$.
\end{theorem}

\section{Preliminary}

\textbf{Notations.}

(a) We use the convention of functional notation, namely, elements of
the mapping class group are applied right to left, i.e. the composition $FG$ means
that $G$ is applied first.

(b) A Dehn twist means a right-hand Dehn twist.

(c) We denote the curves by lower case letters $a$, $b$, $c$, $d$
(possibly with subscripts) and the Dehn twists about them by the corresponding
capital letters $A$, $B$, $C$, $D$. Notationally we do not distinguish a
diffeomorphism/curve and its isotopy class.

\bigskip

\textbf{Cross-cap slide.}

In \cite{Li}, Lickorish proved that the Dehn twists of all the two-sided curves
on the non-orientable surface generate $\mathcal{T}(N_g)$ and $[\text{MCG}(N_g):\mathcal{T}(N_g)] = 2$.
As an example of the mapping classes which do not lie in $\mathcal{T}(N_g)$,
he described a mapping class so-called "Y-homeomorphism" or "cross-cap slide" as shown in figure 1.

\begin{figure}[htbp]
\centering
\includegraphics[totalheight=4cm]{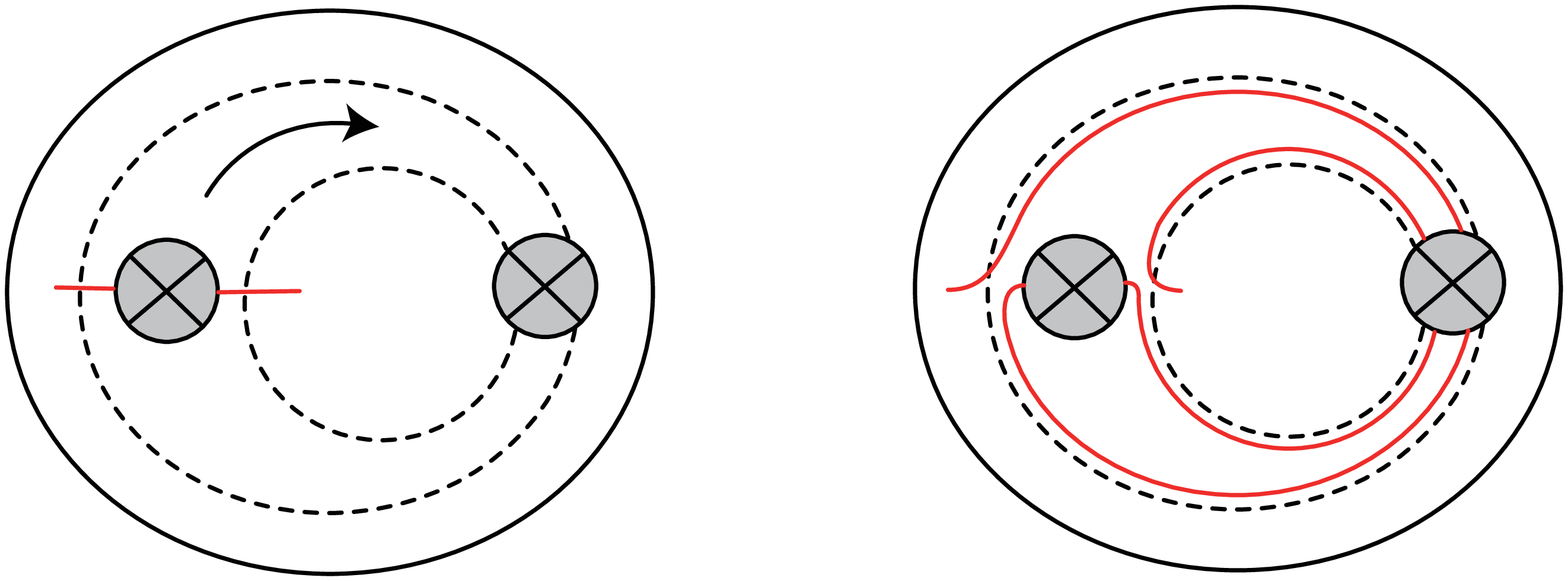}\\
Figure 1.
\end{figure}

\bigskip

\textbf{Two points of view for the M\"obius band partition of a non-orientable surface of odd genus.}

If $g$ is odd, we can decompose the non-orientable surface $N_g$ into $g$ M\"obius bands.
Figure 2 shows two points of view to do this.

(1) The left picture of figure 2 is a $2g$-gon,
with a cross-cap in the middle and the opposite sides glued together pairwise.
Under this gluing, the vertices of this $2g$-gon is divided into two equivalent classes.
After the gluing, they form two points on $N_g$. We denote them as $N$ and $S$.
There are $g$ arcs connecting pairs of antipodal vertices
and passing the cross-cap in the middle of the $2g$-gon.
They divide $N_g$ into $g$ M\"obius bands.
We call this is the $2g$-gon presentation of $N_g$.

\begin{figure}[htbp]
\centering
\includegraphics[totalheight=4.2cm]{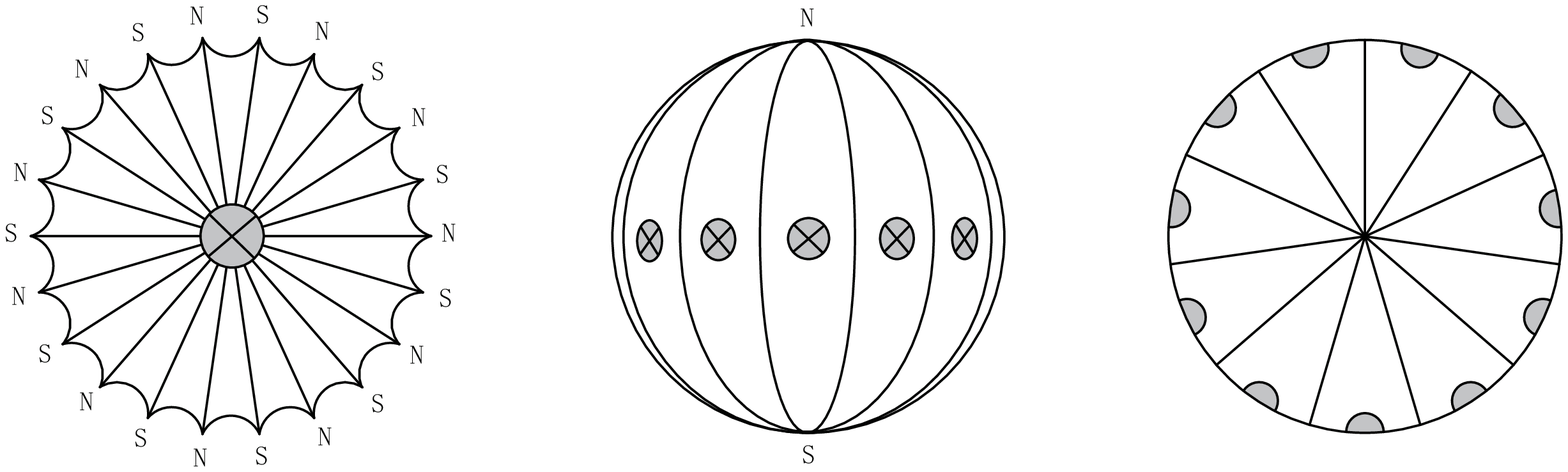} \\
Figure 2.
\end{figure}

(2) The middle and the right picture of figure 2 show a 2-sphere with $g$ projective planes attached.
This is also $N_g$. Suppose the $g$ projective plane sit on the equator.
Denote the north pole and the south pole as $N, S$.
There are $g$ arcs connecting $N$ and $S$. They divide $N_g$ into $g$ M\"obius bands.
We call this is the $g$-cross-cap presentation of $N_g$.

We can check the above two presentations of $N_g$ are equivalent.
In the following, we will go back and forth between such presentations.

\bigskip

\textbf{Mapping classes supported on an Klein bottle with boundary}

For the non-orientable surface $N_g$,
there is a subsurface homeomorphic to a one-holed Klein bottle, see figure 3.
We use the notation as those in \cite{Sz}.
The one-holed Klein bottle contains two cross-caps in its interior.
Suppose $U$ is the mapping class that exchanges the two cross-caps,
which is like a half-twist generator of the braid group.
There is a curve $a$ which is two-sided and passes both cross-caps.
Let $A$ be the Dehn twist along $a$.
The mapping class $Y$ is a $Y$-homeomorphism,
sliding one cross-cap along the one-sided curve which passes the other cross-cap once.

\begin{figure}[htbp]
\centering
\includegraphics[totalheight=5cm]{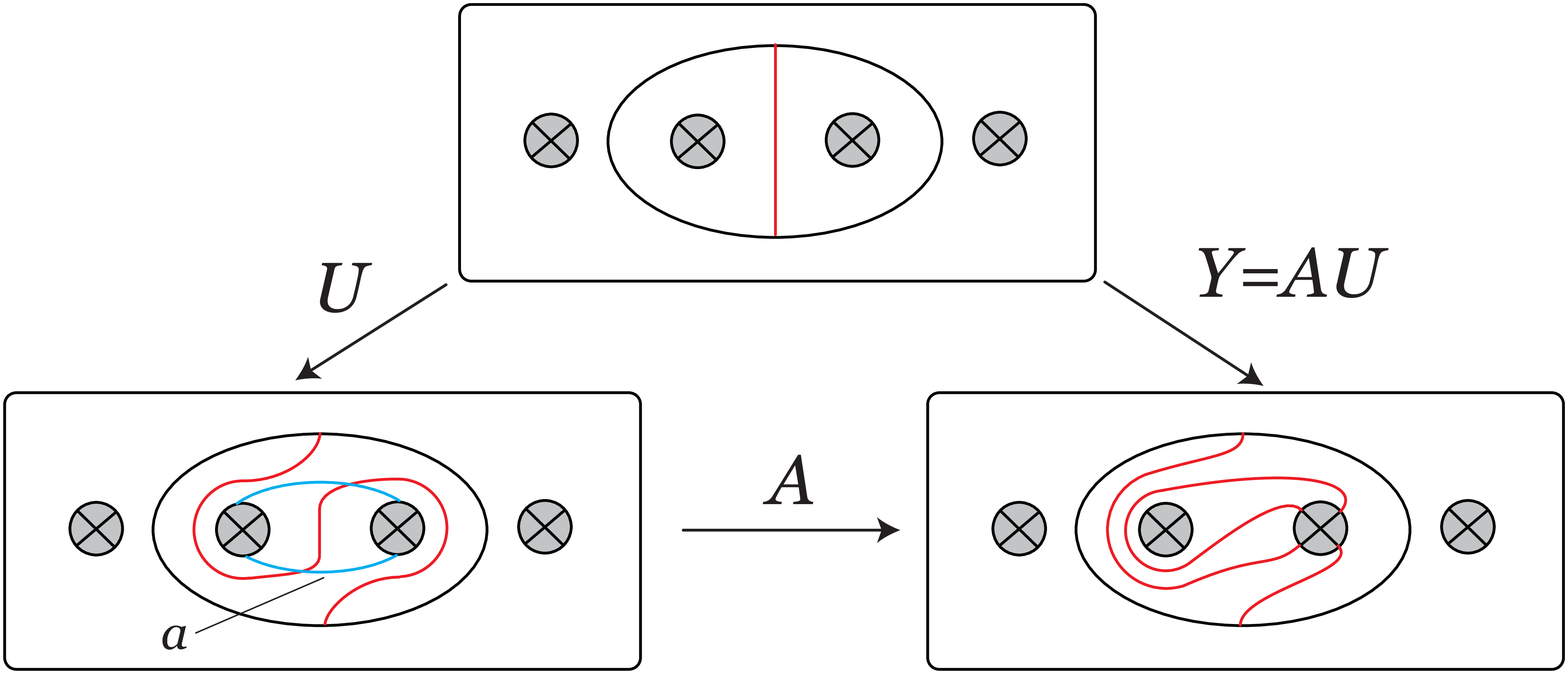}\\
Figure 3.
\end{figure}

Szepietowski showed $Y=AU$. In other words,
both the cross-cap exchanging map $U$ and the cross-cap side $Y$ do not lie in $\mathcal{T}(N_g)$.
Moreover, it is not hard to check that both $Y^2$ and $U^2$
equals the Dehn twist along the boundary curve of the one-holed Klein bottle.

\bigskip

\textbf{The curves needed for generating $\mathcal{T}(N_g)$.}

Omori construct a generating set consist of $g+1$ Dehn twists for $\mathcal{T}(N_g)$ (\cite{Om}).
When we use the $g$-cross-cap presentation of $N_g$,
the curves for those Dehn twists are $a_1, a_2, \dots, a_g, b_0, e$ shown in figure 4.
We can check $A_1^{-1}(e)=c$. Hence the Dehn twists along
$a_1, a_2, \dots, a_g, b_0, c$ can also generate $\mathcal{T}(N_g)$.

\begin{figure}[htbp]
\centering
\includegraphics[totalheight=4cm]{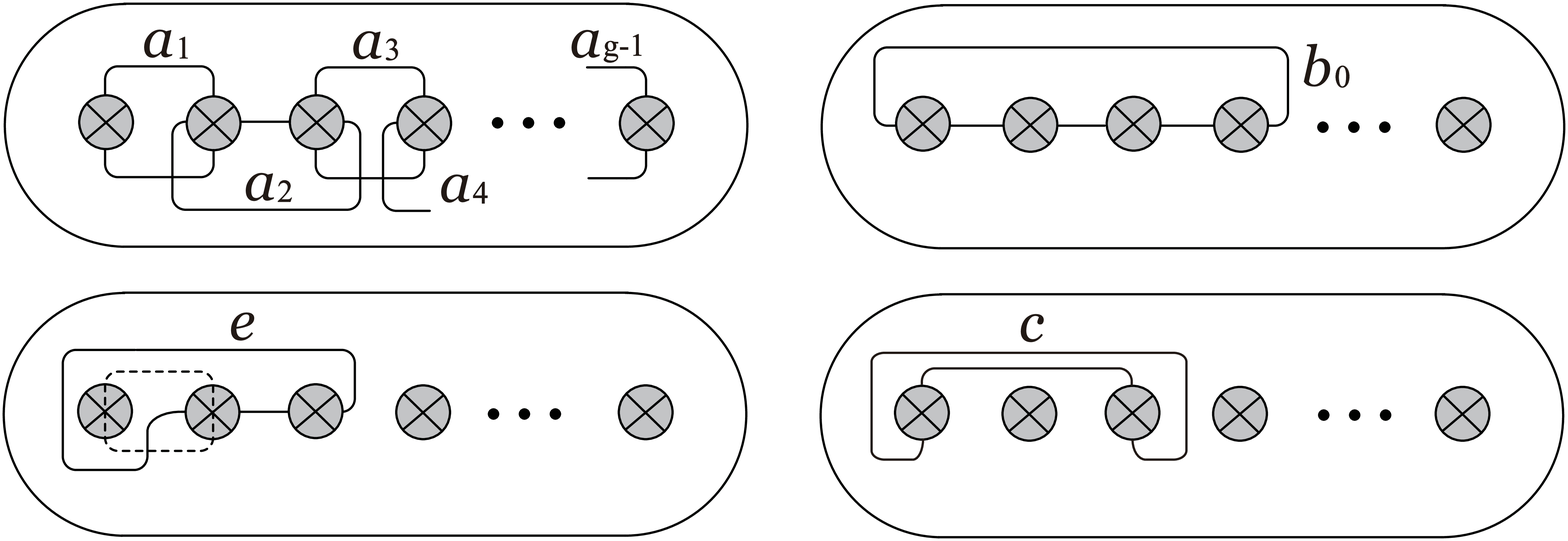} \\
Figure 4.
\end{figure}

We can also use the $2g$-gon presentation to see what these curves are.
See figure 5.

\begin{figure}[htbp]
\centering
\includegraphics[totalheight=3.5cm]{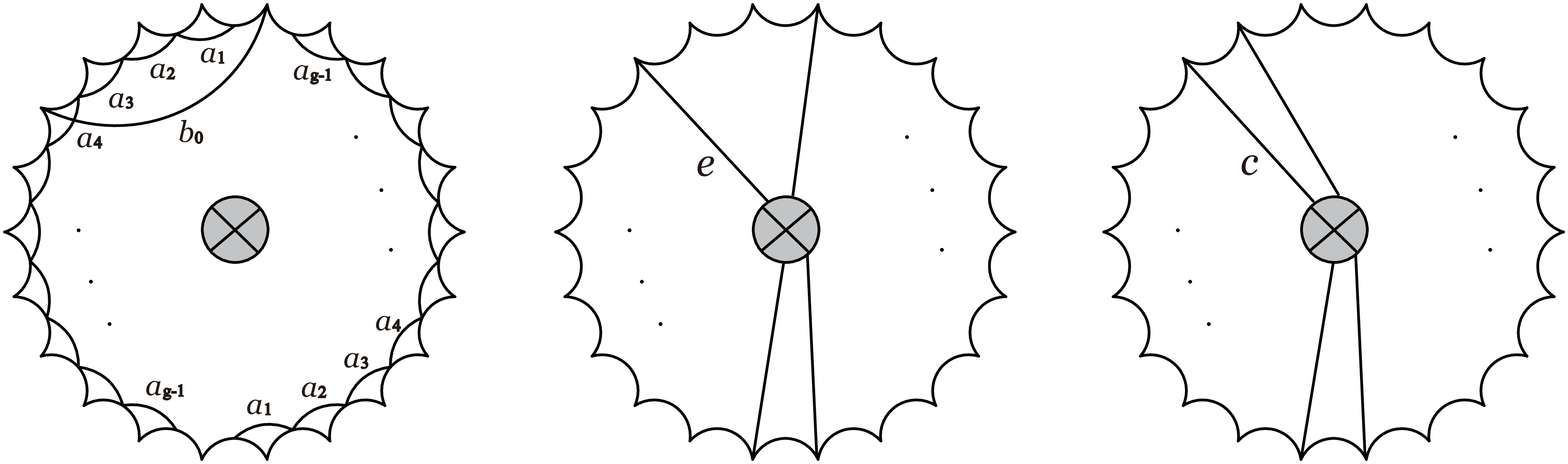} \\
Figure 5.
\end{figure}

\bigskip

\section{The proof of the main theorem}

We now give a proof for Theorem 1.1.

\begin{proof}[Proof of Theorem 1.1]
We first give the torsion generators.
Suppose $g$ is odd. See figure 6.

\begin{figure}[htbp]
\centering
\includegraphics[totalheight=4cm]{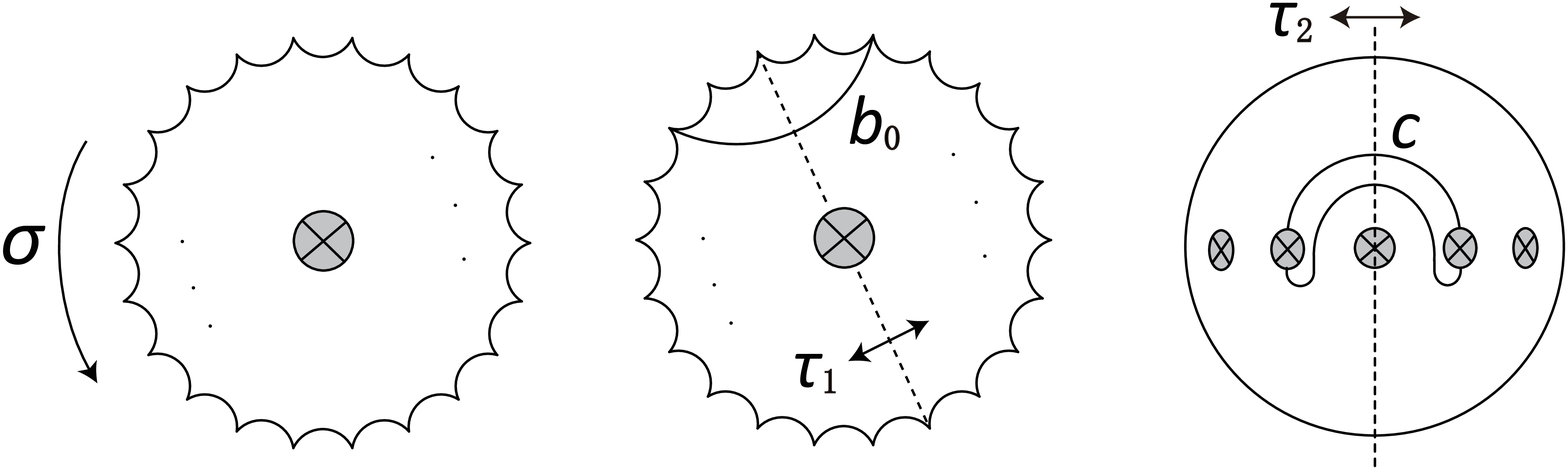} \\
Figure 6.
\end{figure}

Let $\{a_1, \dots, a_g, b_0, c\}$ be the set of curves
whose Dehn twists generate $\mathcal{T}(N_g)$,
$\sigma$ be the rotation in the $2g$-gon presentation,
$\tau_1$ be a reflection of the $2g$-gon presentation that preserves the curve $b_0$.
$\tau_2$ be a reflection of the $g$-cross-cap presentation that preserves $c$.
We can easily see that $\sigma^{2g}=1, (\tau_1 \circ B_0)^2=1, (\tau_2 \circ C)^2=1$.

Let $G = \langle \sigma, \tau_1 \circ B_0, \tau_2 \circ C \rangle$
be the subgroup of $\text{MCG}(N_g)$ generated by these three elements of finite orders.
We claim that when $g$ is odd, $G = \mathcal{T}(N_g)$.

By the method in \cite{Du1} and \cite{Du2},
the Dehn twists $A_1, \dots, A_g, B_0$ are in $G$.
Then $\tau_1$ is also in $G$.

We can interpret some of the torsion elements in more geometric ways. See figure 7.
We can check that $\tau_1$ is not only a reflection in the $2g$-gon presentation
but also a reflection in the $g$-cross-cap presentation.
Let $\tau_3$ be the north-south reflection of the $g$-cross-cap presentation of $N_g$,
$t$ be an order $g$ rotation. Then $\sigma = t \circ \tau_3$ and $\tau_3 = \sigma^g$.
Hence $\tau_3$ and $t$ are also in $G$.

\begin{figure}[htbp]
\centering
\includegraphics[totalheight=3.8cm]{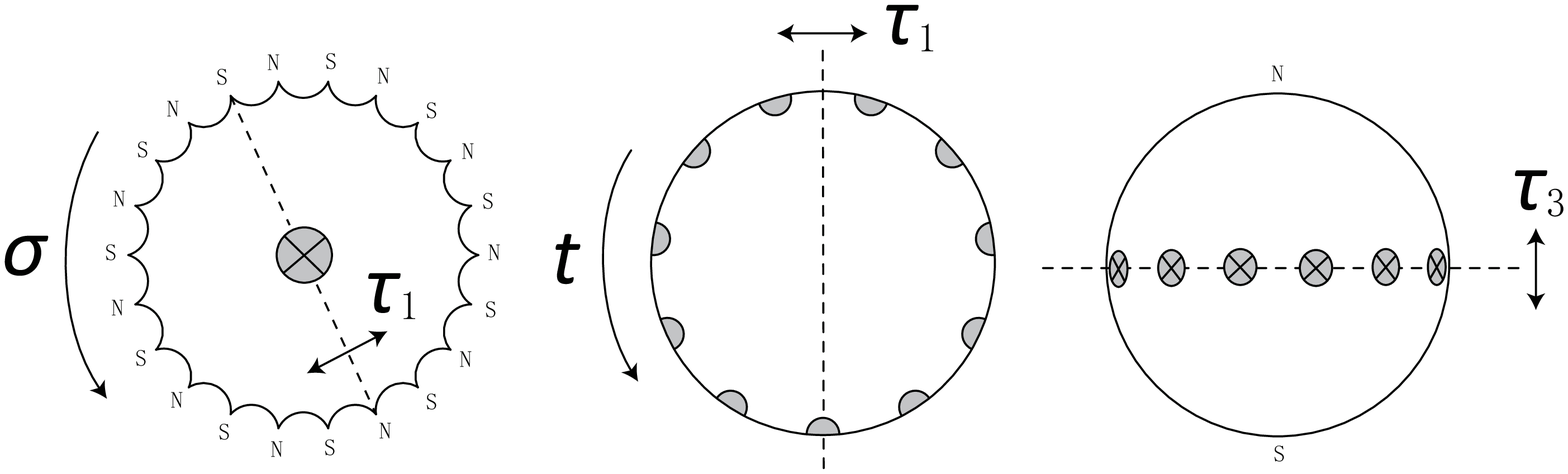} \\
Figure 7.
\end{figure}

Now $\tau_2$ is conjugated to $\tau_1$ by some power of $t$.
So $\tau_2$ also lies in $G$. Hence $C$ lies in $G$.
Since $A_1, \dots, A_g, B_0, C$ generate $\mathcal{T}(N_g)$,
$\text{MCG}(N_g) \geq G \geq \mathcal{T}(N_g)$.
We want to prove $G$ is not $\text{MCG}(N_g)$.
We need to verify all the generators lie in $\mathcal{T}(N_g)$.

Since $\sigma = A_{2g} A_{2g-1} \dots A_2 A_1$, $\sigma$ is in $\mathcal{T}(N_g)$.
So $\tau_3$ and $t$ also lies in $\mathcal{T}(N_g)$.
When $g$ is odd, in the $g$-cross-cap presentation of $N_g$,
the composition of $\tau_1$ and $\tau_3$ is a rotation of the 2-sphere, fixing one cross-cap.
If we look every cross-cap as a punctured point,
then $\tau_1 \circ \tau_3$ becomes an element in the spherical braid group.
It can be written as a product of the standard half-twist generators of the braid group.
We can check the number of the half-twists in the product is an even number.
Each half-twist corresponds to a mapping class exchanging two cross-caps
and supported on a Klein bottle with one boundary.
This means the half-twist generators in the braid group correspond to
the mapping classes outside $\mathcal{T}(N_g)$.
The number of halt-twists in the product is even means
$\tau_1 \circ \tau_3$ lies in $\mathcal{T}(N_g)$, hence $\tau_1$ also lies in $\mathcal{T}(N_g)$.
$\tau_2$ is in $\mathcal{T}(N_g)$ because it is conjugated to $\tau_1$.
We get $G = \langle \sigma, \tau_2 \circ C, \tau_1 \circ B_0 \rangle = \mathcal{T}(N_g)$.
\end{proof}

\end{document}